\documentclass[a4paper, 11pt]{article}
\usepackage{amsfonts}
\usepackage{amsmath}
\usepackage[english]{babel}

\newenvironment{proof}[1][Proof]{\textbf{#1.} }{\ \rule{0.5em}{0.5em}}
\newtheorem{theorem}{Theorem}[section]
\newtheorem{lemma}[theorem]{Lemma}
\newtheorem{e-proposition}[theorem]{Proposition}
\newtheorem{corollary}[theorem]{Corollary}
\newtheorem{e-definition}[theorem]{Definition\rm}

\begin{document}

\title{Well Posedness for Positive Dyadic Model}

\author{D.\ Barbato, F.\ Flandoli and F.\ Morandin}

\maketitle

% 76B03; 35Q31

\begin{abstract}
We consider the solutions of the Cauchy problem for a dyadic model of
Euler equations. We prove global existence and uniqueness of
Leray-Hopf solutions in a rather large class $\mathcal{K}$ that
implies in particular global existence and uniqueness in $l^2$ for all
initial positive conditions in $l^2$.
\end{abstract}

\section{Introduction}
The dyadic model was introduced by Kats and Pavlovic in
\cite{KatPav}. It is a shell-type model that describes the evolution
of wavelet coefficients of Euler and Navier-Stokes solutions. A
construction of the dyadic model can be found in \cite{KatPav},
\cite{Ces} and \cite{Wal}. The infinite system of ODE's of the
inviscid dyadic model is the following:
\begin{equation}
\frac{d}{dt}X_{n}(t)=k_{n-1}X_{n-1}^{2}(t)-k_{n}X_{n}(t)X_{n+1}(t),\qquad
t\geq 0,\qquad X_{n}(0)=x_{n},  \label{eq: dyadic system}
\end{equation}
for $n\geq1$, with $0\leq k_{n}\leq C 2^n$ for every $n\geq
1$, $k_{0}=0$ and $X_{0}(t)=0\ $ for every $t\geq 0$.  Let us
recall or introduce some notations. For the purpose of this paper, the
ambient space for this system is the Hilbert space $H=l^{2}$, whereas
the natural Sobolev space often used in the literature is
$H^{1}=\{u=(u_{n})_{n_{\in }\mathbb{N}}:\Vert u\Vert _{H^{1}}^{2}:=
\sum_{n}k_{n}^{2}u_{n}^{2}<\infty \}$.

\begin{e-definition}
\label{def sol 1} A solution on $[0,T)$ (global solution if $T=+\infty $) of
(\ref{eq: dyadic system}) is a function $X=(X_{n})_{n\in \mathbb{N}}$
defined on $[0,T)$ such that $X_{n}\in C^{1}([0,T))$ for all $n$, $X$
satisfies system (\ref{eq: dyadic system}) and $X(t)\in H$ for all $t\in
\lbrack 0,T)$.
\end{e-definition}

We know that for every initial condition in $H$ there exists at least
one global solution. The uniqueness of the solution in $H^1$ is proved
in \cite{KiZlat}, but it is also shown there that there is a blow-up
in finite time, that is, for all non zero initial condition in $H^1$
there is no global solution in $H^1$.  Similar results are proved in
\cite{KatPav},\cite{FriPav} and \cite{Ces}. On the other hand in
\cite{BFM} there are examples of non-uniqueness in $H$.

Let us restrict the attention to solutions which satisfy suitable
energy inequalities, as initially proposed by \cite{Ler} in the
viscous case. In \cite{DeLellis S 2} it is proved that even in such
classes there are counterexamples to uniqueness.

\begin{e-definition}
\label{def sol 2} A solution $X$ on $[0,T)$, satisfies the \emph{weak energy
inequality} if 
\begin{equation}
\Vert X(0)\Vert _{H}\geq \Vert X(t)\Vert _{H}\qquad \forall t\in \lbrack 0,T),
\label{eq: weak}
\end{equation}
and it satisfies the \emph{strong energy inequality} if 
\begin{equation}
\Vert X(s)\Vert _{H}\geq \Vert X(t)\Vert _{H}\qquad \forall s,t\in [0,T)\ \textup{with}\ s<t.  \label{eq: strong}
\end{equation}
Finally, a solution that satisfies the strong energy inequality will be
called \emph{Leray-Hopf solution}.
\end{e-definition}

\begin{theorem}
\label{th: weak->strong} A solution satisfies the weak energy
inequality if and only if it satysfies the strong energy inequality.
\end{theorem}

This theorem, which is proved at the end of next section, is false for
the original Euler equations. Counterexamples can be found in \cite{DeLellis S 2}.

\begin{theorem}
\label{thm existence_energy_solution} All initial conditions 
in $H$ admit a global Leray-Hopf solution.
\end{theorem}
This theorem is inherited by Navier-Stokes equations \cite{Ler},
\cite{Gal}.  The proof is standard and actually quite simpler in this
setting, since theorem \ref {th: weak->strong} allows to prove just
the weak energy inequality. The details can be found in \cite{BFM}.

\begin{e-definition}\label{def:class_k}
We call \textit{solution of class} $\mathcal{K}$ any solution
$X=(X_{n})_{n\in \mathbb{N}}$ such that the function
\begin{equation}
a( t) =\sup_{n\in \mathbb{N}}( -k_{n}X_{n+1}( t)
)   \label{function a}
\end{equation}
is locally integrable on $[0,\infty )$.
\end{e-definition}

Our main result is the uniqueness of Leray-Hopf solutions of class
$\mathcal{K}$.

\begin{theorem}
\label{th: uniqueness} Let $X^{( i) }=(X_{n}^{( i)
})_{n\in \mathbb{N}}$, $i=1,2$, be two Leray-Hopf solutions with the
same initial condition $x=(x_{n})_{n\in \mathbb{N}}\in H$. Assume that
$X^{( i) }$ are of class $\mathcal{K}$ (in particular this is true if
$x$ has at most a finite number of negative components). Then $X^{( 1)
}=X^{( 2) }$.
\end{theorem}

\begin{corollary}
For every positive initial condition in $H$ there exists a unique
global solution and this is a Leray-Hopf solution.
\end{corollary}

\section{Elementary and known facts}

Denote by $S( t) x$ the set of all values at time $t$ of solutions
with initial condition $x$. This defines a (possibly) multivalued map
$S( t) :H\rightarrow \mathcal{P}( H) $, for all $t\geq 0$, where
$\mathcal{P}( H) $ is the set of all parts of $H$. Let us call $S(t)$
the \textit{multivalued flow} associated to the dyadic model. Let us
introduce the set $H_{+}$ of all $x\in H$ such that $x_{n}<0$ for at
most a finite number of $n$'s. Given a solution $X$, let $E(t)$ be the
energy at time $t$, $E(t) =\sum_{j=1}^{\infty }X_{j}^{2}(t) $ and let
$E_{n}(t) =\sum_{j=1}^{n}X_{j}^{2}( t) $. We begin with a few
elementary properties about solutions and energy.

\begin{lemma}
\label{lemma 1}If $X$ is a solution, then 
\begin{equation*}
X_{j}(s) \geq 0\quad \Longrightarrow \quad X_{j}( t)
\geq 0\text{ for all }t>s.
\end{equation*}
In particular $H_{+}$ is invariant: $S( t) H_{+}\subset H_{+}$.
\end{lemma}

\begin{proof}
By the variation of constants formula, $X_{j}( t) =e^{A_{j}( t) }X_{j}( s)
+\int_{s}^{t}e^{A_{j}( t) -A_{j}( \theta ) }k_{j-1}X_{j-1}^2( \theta )
d\theta $, where $A_{j}( t) =-\int_{s}^{t}k_{j}X_{j+1}( \theta )
d\theta $.
\end{proof}

\begin{lemma}
\label{lemma 2}%
We have $\frac{d}{dt}E_{n}( t)  \leq 0\Leftrightarrow X_{n+1}(
t) \geq 0$ and $
\frac{d}{dt}E_{n}( t)  \geq 0\Leftrightarrow X_{n+1}(
t) \leq 0$.
\end{lemma}

\begin{proof}
It follows from the identity $\frac{d}{dt}E_{n}=-2k_{n}X_{n}^{2}X_{n+1}$.
\end{proof}

\begin{lemma}
\label{lemma 3}If $E( t) <E( s) $ for some $t>s$, then 
$X( t) \in H_{+}$. In plain words: energy may decrease only in
$H_{+}$.
\end{lemma}

\begin{proof}
By contradiction from lemmas \ref{lemma 1} and \ref{lemma 2}: if $X(
t) \notin H_{+}$, there is a sequence $\{n_k\}_k$ such that for every
$k$, $X_{n_k}( t) \leq 0$, yielding $X_{n_k}( \theta ) \leq 0$ for all
$\theta \in [ s,t] $ (lemma \ref{lemma 1}). Hence
$\frac{d}{dt}E_{n_k-1}( \theta ) \geq 0$ for all $\theta \in [ s,t] $
(lemma \ref{lemma 2}), so that $E_{n_k-1}( t) \geq E_{n_k-1}( s)
$. Since $E_n\uparrow E$, this implies $E( t) \geq E( s) $.
\end{proof}

\begin{lemma}
\label{lemma 4}If $E( t) >E( s) $ for some $t>s$, then 
$X( s) \notin H_{+}$. In plain words: energy may increase only in 
$H_{+}^{c}$. In particular, in $H_{+}$ the energy is non-increasing.
\end{lemma}

\begin{proof}
By contradiction from lemmas \ref{lemma 1} and \ref{lemma 2}: if $X(
s) \in H_{+}$, there is an $n_0$ such that for every $n\geq n_0$,
$X_{n}( s) \geq 0$, yielding $X_{n}( \theta ) \geq 0$ for all $\theta
\in [ s,t] $ (lemma \ref{lemma 1}), hence $\frac{d}{dt}E_{n-1}( \theta
) \leq 0$ for all $\theta \in [s,t] $ (lemma \ref{lemma 2}), hence
$E_{n-1}( t) \leq E_{n-1}( s) $, which implies, in the limit as
$n\rightarrow \infty $, $E( t) \leq E( s) $.
\end{proof}

\begin{corollary}
Every solution with initial condition in $H_{+}$ is a Leray-Hopf
solution.
\end{corollary}

\begin{proof}[Proof of theorem \ref{th: weak->strong}]
Let $X$ be a solution, we have to prove that $(\ref{eq:
  weak})\Rightarrow(\ref {eq: strong})$. Let $0<s<t$, from (\ref{eq:
  weak}) we have $E(s)\leq E(0)$ and $E(t)\leq E(0)$. If $E(s)= E(0)$
we are finished.  Otherwise $E(s)< E(0)$, so by lemma \ref{lemma 3},
$X(s)\in H_+$, hence by lemma \ref{lemma 4}, $E(s)\geq E(t)$.
\end{proof}

\section{Uniqueness}
Uniqueness is proved in the class $\mathcal K$ (see
definition~\ref{def:class_k}), which is a technical requirement easily
satisfied if the initial condition is in $H_+$. Truly, if $X(0)\in
H_{+}$, then $X_{n}( t) >0$ for all $n$ larger than some $n_{0}$ and
all $t\geq0$, hence
\begin{equation*}
a( t) \leq \sup_{n\leq n_{0}}k_{n}| X_{n}(
t)|  \leq k_{n_{0}}\sup_{n\leq n_{0}}| X_{n}(
t)| \leq k_{n_{0}}\sqrt{E( t) }\leq k_{n_{0}}\sqrt{E( 0) },
\end{equation*}
where the last inequality is due to lemma \ref{lemma 4}. Hence
solutions starting in $H_{+}$ are of class $\mathcal{K}$. Whether
every Leray-Hopf solution is of class $\mathcal{K}$ is an open
question. We don't know counterexamples and the proof doesn't look
trivial.

It follows from the considerations above, that theorem \ref{th:
  uniqueness} states, in particular, that $S(t)$ is univalued on
$H_{+}$.

\begin{proof}[Proof of theorem \ref{th: uniqueness}]
By hypothesis the energies of $X^{(1)}$ and $X^{(2)}$ are non-increasing
functions of $t$. Hence 
\begin{equation}
| X_{n}^{( i) }(t)| \leq \sqrt{E( 0) },
\label{simple bound}
\end{equation}
for all $n\geq 1$, $t\geq 0$ and $i=1,2$. We shall use this bound below. Let 
\begin{equation*}
Z_{n}:=X_{n}^{( 1) }-X_{n}^{( 2) },\qquad
Y_{n}:=X_{n}^{( 1) }+X_{n}^{( 2) }.
\end{equation*}
It is easy to check that for all $n\geq1$, $Z_{n}(0)=0$ and for $t\geq 0$,
\begin{equation*}
\frac{d}{dt}Z_{n}=k_{n-1}Z_{n-1}Y_{n-1}-\frac{k_{n}}{2}
(Z_{n}Y_{n+1}+Y_{n}Z_{n+1}).
\end{equation*}
This implies 
\begin{equation*}
\frac{d}{dt}
Z_{n}^{2}=2k_{n-1}Y_{n-1}Z_{n-1}Z_{n}-k_{n}Y_{n+1}Z_{n}^{2}-k_{n}Y_{n}Z_{n}Z_{n+1}.
\end{equation*}
One could think of adding up these equations, believing in some
cancellations between two consecutive of them. The terms
$k_{n}Y_{n+1}Z_{n}^{2}$ and $k_{n+1}Y_{n+2}Z_{n+1}^{2}$ have a
dissipative nature. The difficulty in the terms
$-k_{n}Y_{n}Z_{n}Z_{n+1}$ and $2k_{n}Y_{n}Z_{n}Z_{n+1}$ is that they
differ by a factor $-2$. For this reason, instead of using the
classical quantity $\sum_{i=1}^{n}Z_{i}^{2}(t)$, we introduce
\begin{equation*}
\psi _{n}(t):=\sum_{i=1}^{n}\frac{Z_{i}^{2}(t)}{2^{i}}.
\end{equation*}
Indeed
\begin{equation*}
\frac{d}{dt}\frac{Z_{n}^{2}}{2^{n}}=k_{n-1}Y_{n-1}\frac{Z_{n-1}Z_{n}}{2^{n-1}
}-k_{n}Y_{n+1}\frac{Z_{n}^{2}}{2^{n}}-k_{n}Y_{n}\frac{Z_{n}Z_{n+1}}{2^{n}},
\end{equation*}
so $\psi _{n}(t)$ satisfies the simple equation 
\begin{equation*}
\frac{d}{dt}\psi _{n}(t)=-\sum_{i=1}^{n}k_{i}Y_{i+1}\frac{Z_{i}^{2}}{2^{i}}-
\frac{k_{n}}{2^{n}}Y_{n}Z_{n}Z_{n+1},\qquad \psi _{n}(0)=0.
\end{equation*}
Since both solutions are of class $\mathcal{K}$, denoted by $a$ the
maximum of the functions $a^{(i)}$'s defined by (\ref{function a}), we
have $-k_{i}Y_{i+1}( t) \leq 2a( t) $. Hence, using also
$|Z_{n+1}|<2\sqrt{E(0)}$, and $k_n\leq C2^n$, we have
\begin{equation*}
\frac{d}{dt}\psi _{n}(t)
\leq 2a( t) \psi _{n}(t)+KY_nZ_n
\leq 2a( t) \psi _{n}(t)+K\bigl(|X_n^{(1)}|^2+|X_n^{(2)}|^2\bigr).
\end{equation*}
From Gronwall lemma
\begin{equation*}
0\leq\psi _{n}(t)
\leq \int_{0}^{t}e^{\int_{s}^{t}2a( r) dr}K\sum_i|X_n^{(i)}(s)|^2ds
\leq K'(t)\sum_i\int_{0}^{t}|X_n^{(i)}(s)|^2ds,
\end{equation*}
where $K'(t)$ is a positive constant for every $t\geq 0$. Since
$\sum_n\int_{0}^{t}|X_n^{(i)}(s)|^2=\int_{0}^{t}E^{(i)}(s)ds\leq
tE(0)$, both integrals above ($i=1,2$) tend to zero as
$n\rightarrow\infty$, and $\psi _{n}(t)$ as well. Since the latter is
non-decreasing in $n$, we get that $\psi _{n}(t)=0$ for every $n$, and
every $t\geq 0$. This implies $Z\equiv 0$.
\end{proof}
\bigskip

%%\label{}
% etc, etc

% \appendix

% \section{}
% \label{}

% The Acknowledgements are an un-numbered section
%\section*{Acknowledgements}
% Acknowledgements text here


\begin{thebibliography}{00}
% please try to use the bibitem system -
% the references should be in alphabetical order of authors' names.
% Articles with a single author first, author will 1 co-author next,
% then author with several co-authors;


\bibitem{BFM}  D. Barbato, F. Flandoli, F. Morandin, Energy dissipation and
self-similar solutions for an unforced inviscid dyadic model,
arXiv:0811.1689v1.

\bibitem{Ces}  A. Cheskidov, Blow-up in finite time for the dyadic model of
the Navier-Stokes equations, \textit{Trans. Amer. Math. Soc. }\textbf{360}
(2008), no. 10, 5101--5120.

\bibitem{DeLellis S 2}  C. De Lellis, L. Sz\'{e}kelyhidi, On admissibility
criteria for weak solutions of the Euler equations, arXiv:0712.3288.

\bibitem{FriPav}  S. Friedlander, N. Pavlovic, Blowup in a three-dimensional
vector model for the Euler equations, \textit{Comm. Pure Appl. Math.} 
\textbf{57} (2004), no. 6, 705--725.

\bibitem{Gal}  G.P. Galdi, An introduction to the Navier-Stokes
initial-boundary value problem, \textit{Fundamental directions in
mathematical fluid mechanics}, Adv. Math. Fluid Mech. Birkh\"{a}user Basel
(2005), 1-70.

\bibitem{KatPav}  N. H. Katz, N. Pavlovic, Finite time blow-up for a dyadic
model of the Euler equations, \textit{Trans. Amer. Math. Soc.} \textbf{357}
(2005), no. 2, 695--708.

\bibitem{KiZlat}  A. Kiselev, A. Zlato\v{s}, On discrete models of the Euler
equation, \textit{IMRN} \textbf{38} (2005), no. 38, 2315-2339.

\bibitem{Ler}  J. Leray, Sur le mouvement d'un liquide visqueux emplissant
l'espace, \textit{Acta Math.} \textbf{63} (1934), no.1 193-248.

\bibitem{Wal}  F. Waleffe, On some dyadic models of the Euler equations, 
\textit{Proc. Amer. Math. Soc.} \textbf{134} (2006), 2913-2922.
\end{thebibliography}
\end{document}